\documentclass{article}

\usepackage{stywhispers,epsfig}
\usepackage[named]{algo}
\usepackage{amsmath}
\usepackage{amsfonts}
\usepackage{amssymb}
\usepackage{graphicx}
\usepackage{subfigure}
\usepackage{amsthm}
\usepackage{url}

\ninept

\newcommand{\bm}[1]{\mbox{\boldmath $#1$}}

\title{A Variable Splitting  Augmented Lagrangian approach to Linear Spectral Unmixing}
\name{Jos\'{e} M. Bioucas-Dias\thanks{This work was supported by the European Commission
Marie Curie training grant MEST-CT-2005-021175. Email:{\{bioucas\}@lx.it.pt}}}
\address{Instituto de Telecomunica\c{c}\~{o}es,\\
Instituto Superior T\'{e}cnico, Technical University of Lisbon,\\
 Lisboa, Portugal}
\begin{document}
%
\maketitle
\begin{abstract}

This paper presents a new  linear  hyperspectral unmixing method of the minimum volume  class,
termed \emph{simplex identification via split augmented Lagrangian} (SISAL). 
Following  Craig's seminal ideas,
hyperspectral linear unmixing amounts to finding the minimum volume
simplex containing the hyperspectral vectors. This is a nonconvex optimization problem
with convex constraints. In the proposed approach,  the positivity constraints, forcing  the spectral vectors to belong to the convex hull
of the endmember signatures, are replaced by  soft constraints.  The obtained
problem is solved by a sequence of augmented Lagrangian optimizations. The resulting algorithm is very fast and able so solve problems far beyond the reach of the current state-of-the art algorithms. The effectiveness of SISAL is illustrated with
simulated data.

\end{abstract}
\section{Introduction}
\label{sec:intro}

Hyperspectral unmixing is a source separation  problem \cite{ke00}. Compared
with the canonical source separation scenario, the sources in hyperspectral
unmixing (\emph{i.e.}, the materials present in the scene) are statistically
dependent and combine in a linear or nonlinear fashion. These characteristics,
together with the high dimensionality of hyperspectral vectors, place the
unmixing of hyperspectral mixtures beyond the reach of most source separation
algorithms, thus fostering  active research in the field \cite{JNBD05}.

Given a set of mixed hyperspectral vectors, linear mixture analysis, or linear
unmixing, aims at estimating the number of reference materials, also called
endmembers, their spectral signatures, and their abundance fractions
\cite{ke00, JNBD05, plaza04, JN05, LM07, MC94}. The approaches to hyperspectral
linear unmixing can be classified as either  statistical  or  geometrical. The
former  address spectral unmixing as an inference problem, often formulated
under the Bayesian framework, whereas the latter exploit the fact that the
spectral vectors, under the linear mixing model, are in a simplex whose
vertices represent the sought endmembers.

\subsection{Statistical approach to spectral unmixing}

Modeling the  abundance fractions  (sources) statistical dependence  in
hyperspectral unmixing is a central issue in the statistical framework. In
\cite{Dias}, the abundance fractions are modeled as mixtures of Dirichlet densities.
The resulting algorithm, termed DECA (dependent component analysis),
implements an expectation maximization iterative scheme for the inference of
the endmember signatures (mixing matrix) and the  density parameters of the
abundance fractions.

The inference engine in the Bayesian framework is the  posterior density of the
entities to be estimated, given the observations. According to Bayes'  law,
the posterior includes two factors: the  observation density, which may account
for additive noise, and a prior, which may impose constraints on the endmember
matrix ({\em e.g.,} nonnegativity of its elements) and on the abundance
fractions ({\em e.g.,} to be in the probability simplex) and model   spectral
variability. Works \cite{Dobigeon:ieeesp:08, Moussaoui:nc:08} are
representative of this line of attack.

\subsection{Geometrical approach to spectral unmixing}

The geometrical approach exploits the fact that, under the linear mixing model,
hyperspectral vectors belong to a simplex set whose vertices correspond to the
endmembers. Therefore, finding the endmembers is equivalent to identifying the
vertices of the referred to simplex.

If  there exists at least one pure
pixel per endmember (\emph{i.e.}, containing just one material), then unmixing amounts to finding the spectral vectors in the data set corresponding to the vertices of the data simplex. Some popular algorithms taking this assumption are the the N-FINDR \cite{ME99}, the the
\emph{pixel purity  index} (PPI) \cite{JB93}, the \emph{automated morphological
endmember extraction} (AMEE) \cite{plaza02}, the \emph{vertex component
analysis} (VCA) \cite{JN05}, and   the \emph{simplex growing algorithm} (SGA)
\cite{CHANG06}.

If the pure pixel assumption is not fulfilled, what is a more realistic scenario, the
unmixing process is a rather challenging task, since the endmembers, or at least some of
them, are not in the data set. A possible line of attack, in the vein of the seminal
ideas introduced in \cite{MC94}, is to fit a simplex of minimum volume to the data set.
Relevant works exploiting this direction are the \emph{minimum volume enclosing simplex}
(MVES) \cite{Chan_et_al_2009}, the \emph{minimum volume simplex analysys} (MVSA) \cite{li_bioucas08}, and the \emph{nonnegative matrix factorization minimum volume transform} (NMF-MVT) \cite{XT07}.
MVES and MVSA, although implemented in rather different ways, yield state-of-the-art
results. Their major shortcoming is the time they take for more than, say, 10 endmembers
and more than 5000 spectral vectors.

\subsection{Proposed approach}

We introduce the  \emph{simplex identification via split augmented Lagrangian} (SISAL) algorithm for unsupervised hyperspectral linear unmixing.  SISAL belongs to the
minimum volume class, and thus is able to  unmix hyperspectral data sets in which the
pure pixel assumption is violated.

In SISAL, the  positivity  hard constraints are replaced by hinge type soft constraints,
whose strength is controlled by a regularization parameter. This replacement has three advantages: 1)  robustness to outliers and noise; 2) robustness  to poor
initialization; 3) opens the door to dealing  with large problems.  Furthermore, for large values of the regularization parameter, the hard constraint formulation is recovered.

To tackle the hard nonconvex optimization problem we have in hands, we solve a sequence
of nonsmooth convex   subproblems, using  variable splitting to obtain a constraint formulation,
and then applying  an  augmented Lagrangian technique \cite{MarioBioucasManyaTR09,GoldsteinOsher}.  This sub-problems implement an alternate minimization scheme, with  very simple and fast steps.

The paper is organized as follows. Section 2 formulates the problem.
Section 3 proposes a sequence of nonsmooth subproblems, Section 4 solves  the nonsmooth problems
via a variable splitting augmented Lagrangian scheme, Section 5 presents results, and Section 6 ends the paper by presenting a few concluding remarks.

\section{Problem formulation}
\label{sec:MVSA}

Let us assume that in a given scene there are $p$ materials, termed endmembers,   with
spectral signatures ${\bf m}_i\in\mathbb{R}^l$, for $i=1\dots,p$, where $l\geq p$ denotes
the number of spectral bands.  Under the linear mixing model, a given   hyperspectral observed
vector  is a linear combination of the endmember spectral signatures,
where the weights  represent the fractions that each material occupies in the pixel.
Therefore, the observed spectral vectors are in the convex hull of endmember spectral
signatures.

To fix notation, let ${\bf Y}\equiv[{\bf y_{1}},\dots, {\bf y}_{n}]\in\mathbb{R}^{l\times
n}$ denote  a matrix holding  the observed  spectral vectors  ${\bf
y}_{i}\in\mathbb{R}^{l}$ and ${\bf S}\equiv[{\bf s}_1,\dots,{\bf
s}_n]\in\mathbb{R}^{p\times n}$ a  matrix holding the respective
fractions; {\em i.e.},  ${\bf y}_i={\bf M}{\bm s}_i$, for $i=1,\dots,n$, where ${\bf
M}\equiv[{\bf m}_{1},\ldots,{\bf m}_{p}]\in\mathbb{R}^{p\times p}$ is the mixing matrix
containing the endmembers, and ${\bm s}_i$  is a vector denoting  the fractions, often
termed fractional abundances. Since the components of ${\bf s}_i$ are nonnegative and sum
one (they are fractions), then the  fractional abundance vectors belong to the  standard
$p$-simplex set ${\cal S}_p=\{{\bf s}\in\mathbb{R}^p\,:\,{\bf s} \succeq {\bf 0},\,{\bf
1}^T_p{\bf s}=1\}$\footnote{${\bf x}\succeq{\bf y}$ means  $x_i\geq y_i$
for $i=1,\dots,p$; \  ${\bf 1}^T_p\equiv(1,\dots,p)$. }. Therefore,
\begin{equation}\label{obs_model} 
   \begin{array}{lclll}
          {\bf Y} &  = &  {\bf MS}, &  {\bf S} \in{\cal S}_p^n.
   \end{array}
\end{equation}

\begin{figure}[hbt]
    \centering
    \includegraphics[width=40mm,angle=0]{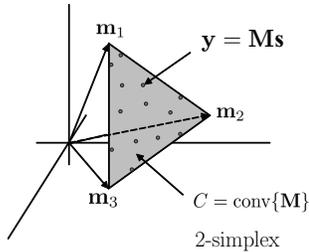}
   \caption{Illustration of the 2-simplex set generated the columns of $\bf M$.}
  \label{fig:convex_set}
\end{figure}

Assuming that the endmember spectral signatures ${\bf m}_i$, for $i=1,\dots,p$, are linearly
independent, then the set $ \{{\bf y}\in\mathbb{R}^l\,:\,{\bf y}={\bf Ms},\;{\bf
s}\in{\cal S}_p\}$ is a $(p-1)$-dimensional simplex, and estimating $\bf M$ amounts to
infer its vertices. Figure \ref{fig:convex_set} illustrates this  perspective.

In this work, we assume that the  number of endmembers and  signal subspace is known
before hand (see, {\em e.g.}, \cite{HYSIME08})  and that the observed vectors ${\bf
y}_i$, $i=1,\dots,n$, represent  coordinates  with respect to a $p$-dimensional basis of
the signal subspace, {\em i.e.},  $l=p$.

Given $\bf Y$, and inspired by the seminal  work \cite{MC94}, we infer matrices $\bf M$
and $\bf S$ by fitting a minimum volume simplex to the data subject to the constraints
${\bf S}\succeq 0$ and ${\bm 1}^T_p{\bf S}={\bf 1}_n$. Since the volume defined by the
columns of $\bf M$ is proportional to $|\mbox{det}\,({\bf M})|$ (note that we are assuming that $\bf M$ is square), then
\begin{equation}\label{opt} 
      \begin{array}{rl}
        {\bf M}^{*} =    & \displaystyle\arg\min_{\bf M} |\det({\bf M})|\\
\mbox{s.t.\ :} & {\bf QY}\succeq {\bf 0},\;\;1_{p}^{T}{\bf QY}={\bf 1}_{n}^{T},
    \end{array}
\end{equation}
where ${\bf Q}\equiv {\bf M}^{-1}$. Since $\det({\bf Q})= 1/\det({\bf M})$,  we can replace the problem
(\ref{opt}) with
\begin{equation}\label{opt2} 
      \begin{array}{rl}
        {\bf Q}^{*}   =    &\displaystyle \arg\min_{\bf Q}  \; -\log |\det({\bf Q})|\\
\mbox{s.t.:}  & {\bf QY}\succeq  {\bf 0},\;\;{\bf 1}_{p}^{T}{\bf QY}={\bf 1}_{n}^{T}.
    \end{array}
\end{equation}

The constraints in (\ref{opt2}) define a convex set. If matrix $\bf Q$ is symmetric and
positive-definite , the problem  (\ref{opt2}) is convex. However, in most cases of
practical interest $\bf Q$ is neither  symmetric nor  positive-definite and, thus, the (\ref{opt2}) is nonconvex.   Therefore, there is no hope in finding systematically the global optima of
(\ref{opt2}). The SISAL algorithm, we introduce below,   aims at ``good'' sub-optimal
solutions.

Our first step  is to simplify the set of constraint ${\bf 1}_{p}^{T}{\bf QY}={\bf
1}_{n}^{T}$. We note that  the vector ${\bf 1}_{n}^T$ does not belong to the null space
of $\bf Y$. Otherwise, the null vector would belong to the the affine hull of $\bf M$,
what would imply that the columns of $\bf M$ would not be independent. Therefore, by
multiplying the equality constraint on the right hand side by ${\bf Y}^T({\bf Y}{\bf
Y}^T)^{-1}$, we get $({\bf 1}_{p}^{T}{\bf QY} ={\bf 1}_{n}^{T}) \Leftrightarrow ({\bf
1}_{p}^{T}{\bf Q}={\bf a}^T)$, where ${\bf a}^T\equiv {\bf 1}_n{\bf Y}^T({\bf Y}{\bf
Y}^T)^{-1}$. The problem (\ref{opt2}) simplifies, then, to
\begin{equation}\label{modprob2} 
      \begin{array}{rl}
      {\bf Q}^{*} =    &\displaystyle\ \arg\min_{\bf Q}\; -\log |\det({\bf Q})|\\
\mbox{s.t.\ :}  & {\bf QY}\succeq {\bf 0},\;\;{\bf 1}_{p}^{T}{\bf Q}={{\bf a}^T}.
    \end{array}
\end{equation}
Instead of solving (\ref{modprob2}), we solve the following modified version:
\begin{equation}\label{modprob_soft} 
      \begin{array}{rl}
      {\bf Q}^{*} =    &\displaystyle\ \arg\min_{\bf Q}\; -\log |\det({\bf Q})|+\lambda \|{\bf QY}\|_h\\
       \mbox{s.t.\ :}  & {\bf 1}_{p}^{T}{\bf Q}={{\bf a}^T},
    \end{array}
\end{equation}
where $\|{\bf X}\|_h \equiv \sum_{ij}h([{\bf X}]_{ij})$ and $h(x)\equiv\max\{-x,0\}$ is the
so-called hinge function. Notice that $\|{\bf QY}\|_h$ penalizes the negative components
of $\bf QY$ proportionally to their magnitude, thus playing the rule of a soft constraint or a {\em regularizer}.  The  amount of regularization is controlled by the
\emph{regularization parameter} $\lambda >0$.
As already referred to, the soft constrained formulation  yields solutions that are
robust to outliers, noise, and poor initialization. Furthermore, by replacing $n\times p$ equality
constraints by a {\em regularizer}, it opens the door to deal with large scale problems.

\section{Sequence of convex subproblems}

Let ${\bf q}\equiv \mbox{vec}({\bf Q})$ denote the operator that stacks the columns of $\bf Q$
in the column vector $\bf q$. Given that
$\mbox{vec}\,({\bf AB})= ({\bf B}^T\otimes {\bf I})\,\mbox{vec}\,({\bf A})= ({\bf
I}\otimes{\bf A})\,\mbox{vec}({\bf B})$, where $\otimes$ denotes the kronecker operator,
and defining $f({\bf q})=-\log |\det({\bf Q})|$,
then (\ref{modprob_soft}) may be written as
\begin{equation}\label{modprob_soft_q} 
      \begin{array}{rl}
      {\bf q}^{*} =    &\displaystyle\ \arg\min_{\bf q}\ f({\bf q})+\lambda \|{\bf Aq}\|_h\\
       \mbox{s.t.\ :}  & {\bf B}{\bf q}={{\bf a}},
    \end{array}
\end{equation}
where ${\bf A}= ({\bf Y}^T\otimes{\bf I})$ and ${\bf B} = ({\bf I}\otimes {{\bf
1}^T_p})$. The Hessian of $f$ is ${\bf H}={\bf K}_n[{\bf Q}^{-T}\otimes {\bf Q}^{-1}]$,
where ${\bf K}_n$ id the comutation matrix ({\em i.e.}, ${\bf K}_n\mbox{vec}({\bf A}) =
\mbox{vec}({\bf A}^T)$). Since $\bf H$  has positive and negative eigenvalues, the
above problem in nonconvex and thus hard.

Using  a  quadratic approximation for $f({\bf q})$,
we approximate (\ref{modprob_soft_q}) by  computing  a descent sequence ${\bf q}_k$, $k=0,1,\,\dots,$ with the following algorithm:
\begin{algorithm}{1 Sequence of strictly convex subproblems}{
\label{alg:alg1}}
Set $k=0$, choose $\mu > 0$ and ${\bf q}_0 =\,\mbox{VCA}(\bf Y)$.\\
\qrepeat\\
     $l_k=f({\bf q}_{k}) + \lambda\|{\bf  Aq}_{k}\|_h$\\
     ${\bf g} = -\mbox{vec}\,({\bf Q}^{-1})$\\
     ${\bf q}_{k+1}  \in \arg\min_{{\bf q}} {\bf g}^T{\bf q}+\mu\|{\bf q}-{\bf q}_k\|^2 +\lambda\|{\bf
     Aq}\|_h$\\
     \hspace{1.0cm}{s.t.:}$\;\;{\bf Bq}={\bf a}$\\
     \qif $f({\bf q}_{k+1}) + \lambda\|{\bf Aq}_{k+1}\|_h > l_k$\\
           $\mbox{find}\; {\bf q}\in \{\alpha {\bf q}_{k+1}+(1-\alpha) {\bf q}_k\,:\,0<\alpha< 1 \}$\\
           \hspace{0cm} $\mbox{such that}\;\; f({\bf q}) + \lambda\|{\bf  Aq}\|_h\leq
           l_k$\\
           ${\bf q}_{k+1}={\bf q}$\qfi\\
     $k \leftarrow k + 1$
\quntil a stopping criterion is satisfied.
\end{algorithm}
\vspace{0.3cm}

Algorithm 1 is initialized with the   VCA \cite{JN05} estimate.
Line 4 computes the gradient of $f({\bf q})$. Line 5 minimizes a
strictly convex approximation to the initial objective function, where the term $f({\bf q})$ was replaced by a quadratic approximation.  The term $\|{\bf q}-{\bf q}_k\|^2$ ensures that $\|{\bf q}_{k+1}-{\bf q}_k\|^2$ does not grow unbounded. Lines 7 to 10
ensures that the objective function  does not increase. To
solve the minimization 5-6, we  introduce in the next  subsection a variable splitting augmented Lagrangian algorithm.

\section{Variable splitting and augmented Lagrangian}

The optimization problem 5-6 of Algorithm 1  is equivalent to
\begin{eqnarray}\label{var_split} 
           \min_{{\bf q}, {\bf z}} &&E({\bf q},{\bf z})\\
       \mbox{s.t.\ :}  & &{\bf Bq}={\bf a},\;\;{\bf Aq}={\bf z}, \nonumber
\end{eqnarray}
where
\[
  E({\bf q},{\bf z}) \equiv {\bf g}^T{\bf q}+\mu\|{\bf q}-{\bf q}_k\|^2
            +\lambda\|{\bf z}\|_h.
\]
In (\ref{var_split}), the variable $\bf q$ was split  into  the pair $({\bf q}, {\bf
z})$ and linked through the constraint ${\bf Aq}= {\bf z}$. The so-called augmented Lagrangian (AL) for this problem, with respect to the constraint ${\bf Aq}={\bf z}$ is given by
\begin{eqnarray}\label{AL}
    {\cal L}({\bf q},{\bf z},{\bf d},\tau) & \equiv& E({\bf q},{\bf z}) + \bm{\alpha}^T({\bf Aq}-{\bf z})+\tau\|{\bf Aq}-{\bf z} \|^2\\
       & = & E({\bf q},{\bf z})  + \tau\|{\bf Aq}-{\bf z-{\bf d}} \|^2 + c,
\end{eqnarray}
where  $\bm{\alpha}$ is holds the Lagrange multipliers, ${\bf d}=-\bm{\alpha}/(2\tau)$, and $c$ is an irrelevant constant. The AL
algorithm consists in minimizing $\cal L$ with respect to $({\bf q}, {\bf z})$ and then updating $\bm{\alpha}$, or, equivalently $\bf d$, as follows
\begin{algorithm}{ 2 Augmented Lagrangian Algorithm}{
\label{alg:alg2}}
Set $t=0$, choose  $({\bf q}_0, {\bf z}_0)$, $\bm{\alpha}_0$, and $\tau > 0$,\\
\qrepeat\\
     $({\bf q}_{t+1}, {\bf z}_{t+1}) \in \arg\min_{{\bf z}} {\cal L}({\bf q},{\bf z},{\bf d}_t,\tau)$\\
     \hspace{1.5cm}s.t.:\;\;${\bf Bq}={\bf a}$\\
     ${\bf d}_{t+1} \leftarrow {\bf d}_{t} - ({\bf Aq}_{t+1} - {\bf z}_{t+1})$\\
     $t \leftarrow t + 1$
\quntil a stopping criterion is satisfied.
\end{algorithm}
\vspace{0.3cm}

It has been shown that, with adequate initializations, the AL
algorithm generates the same sequence as a {\it proximal point algorithm}
(PPA) applied to the Lagrange dual of problem (\ref{var_split});
for further details, see \cite{Iusem, Setzer} and references therein.
Moreover,  the sequence ${\bf d}_{k}$, for $k=0,1,\dots$, converges
to a solution of this dual problem and  all cluster points of
the sequence ${\bf z}_{k}$, for $k=0,1,\dots$, are solutions of the (primal) problem
(\ref{var_split}) \cite{Iusem}.

The exact solution of the optimization with respect to $({\bf q}, {\bf z})$
in the line 3 of the Algorithm 2 is stil a complex task. However, the
block minimizations  with respect to $\bf q$ and  with respect to $\bf z$  are very light
to compute. Based on this,  we propose the following modification of Algorithm 2:
\begin{algorithm}{ 3 Alternating Split AL}{
\label{alg:salsa1}}
Set $t=0$, choose  $({\bf q}_0, {\bf z}_0)$, $\bm{\alpha}_0$, and $\tau > 0$,\\
\qrepeat\\
     $\displaystyle {\bf q}_{t+1} \in \arg\min_{\bf q} {\bf g}^T{\bf q}+\frac{\mu}{2}\|{\bf q}-{\bf q}_k\|^2+ \frac{\tau}{2}\|{\bf Aq}-{\bf z}_t-{\bf d}_t \|^2$\\
     \hspace{1.5cm}s.t.:\;\;${\bf Bq}={\bf a}$\\
     $\displaystyle {\bf z}_{t+1} \in \arg\min_{\bf z} \frac{1}{2} \|{\bf Aq}_{t+1}-{\bf z}-{\bf d}_t \|^2+\frac{\lambda}{\tau}\|{\bf z}\|_h$\\
     ${\bf d}_{t+1} \leftarrow {\bf d}_{t} - ({\bf Aq}_{t+1} - {\bf z}_{t+1})$\\
     $t \leftarrow t + 1$
\quntil stopping criterion is satisfied.
\end{algorithm}
\vspace{0.3cm}
The solution of the  quadratic problem with linear constraints 3-4 is
\begin{equation}
  \label{eq:quadractic_step}
  {\bf q}_{t+1}={\bf F}^{-1}{\bf b}-{\bf F}^{-1}{\bf B}^T({\bf B}{\bf F}^{-1}{\bf B}^T)^{-1}({\bf B}{\bf F}^{-1}{\bf b}-{\bf a}),
\end{equation}
where
\begin{equation}
\begin{array}{lcl}
  {\bf F} & \equiv & (\mu{\bf I}+\tau {\bf A}^T{\bf A}) \\
  {\bf b} & \equiv & \mu {\bf q}_t-{\bf g}+\tau{\bf A}^T({\bf z}_t+{\bf d}_t).
\end{array}
\end{equation}

The minimization with respect to ${\bf z}$, in line 5, is, by definition the \emph{proximity operator of} of the convex function $\|{\bf z}\|_h$  \cite{CombettesSIAM}, which is similar to the soft threshold function but
applied just to the negative part of its argument:
\begin{equation}
\label{eq:IST_step}
{\bf z}_{t+1} = \mbox{soft}_{-}({\bf A}{\bf q}_{t+1}-{\bf d}_t,\mu/\tau),
\end{equation}
where $\displaystyle\mbox{soft}_{-}(x,\beta)=(\max\{|x+\beta/2|-\beta/2,0\})(x/|x|)$ is
a thresholding/shrinkage function and, for a matrix $\bf X$,  $\displaystyle\mbox{soft}_{\_\_}({\bf X},\beta)$ is the  componentwise application of $\displaystyle\mbox{soft}_{-}(\cdot,\beta)$.
We note that computations (\ref{eq:quadractic_step}) and (\ref{eq:IST_step}) are very light: In the first case, all matrices involved are $p^2\times p^2$ and are not iteration dependent. In the second case, the function $\mbox{soft}_{-}(\cdot)$ takes negligible time.

The iterations of  Algorithm 3 are much faster than that of Algorithm 2. There is, however,
the question of convergence. The answer to this question
turns out to be positive, a result that can be proved via the
equivalence between the {\it alternating split AL algorithm} just
described and the so-called Douglas-Rachford splitting method,
applied to the dual of  problem (\ref{var_split}); see \cite[Theorem 8]{EcksteinBertsekas}, \cite{Setzer}.

In conclusion,  the pseudo-code for  SISAL  is given by Algorithm 1 with the step 5 replaced by
Algorithm 3.

\section{Experimental Results}
\label{sec:simu}

This section presents results obtained with SISAL, MVSA (hard version) \cite{li_bioucas08},   MVES \cite{Chan_et_al_2009}, and VCA \cite{JN05}  applied to simulated data sets.  The SISAL regularization parameters was set to  $\lambda=10$,  The remaining parameters were set to $\tau = 1$ and $\mu = 10^{-4}$.
Although these values may be far from optimal, they led to excellent results.
The data was generated according to the linear observation model (\ref{obs_model}).
The abundance fractions are Dirichlet distributed
with  parameter $\mu_{i}=1$, for $i = 1, \ldots, p$. The mixing matrix $\bf M$ is randomly
generated with i.i.d. uniformly distributed elements.
To ensure that no pure pixel is present, we discarded all pixels with any abundance
fractions larger than 0.8. The signal-to-noise ratio  (SNR) defined as $\|{\bf Y}\|^2_F/\|{\bf N}\|_F^2$, where $\|\cdot\|^2_F$ denotes the Frobenius norm and $\bf N$ is zero-mean Gaussian additive noise, was set to SNR=$40$dB.

\begin{figure}[hbt]
    \centering
     \subfigure[]{
    \includegraphics[width=60mm,angle=0]{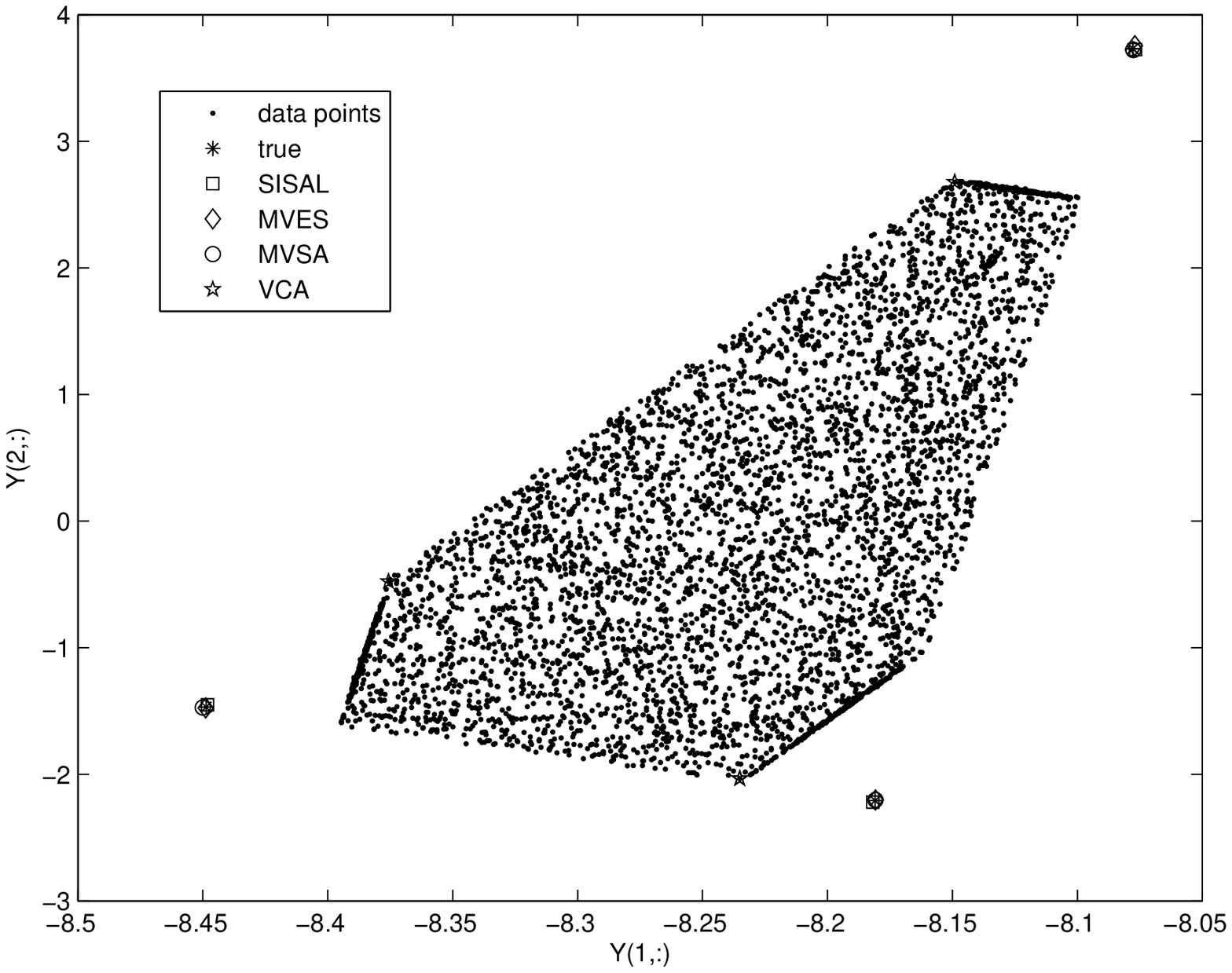}
    }
    \hfill
    \subfigure[]{
    \includegraphics[width=60mm,angle=0]{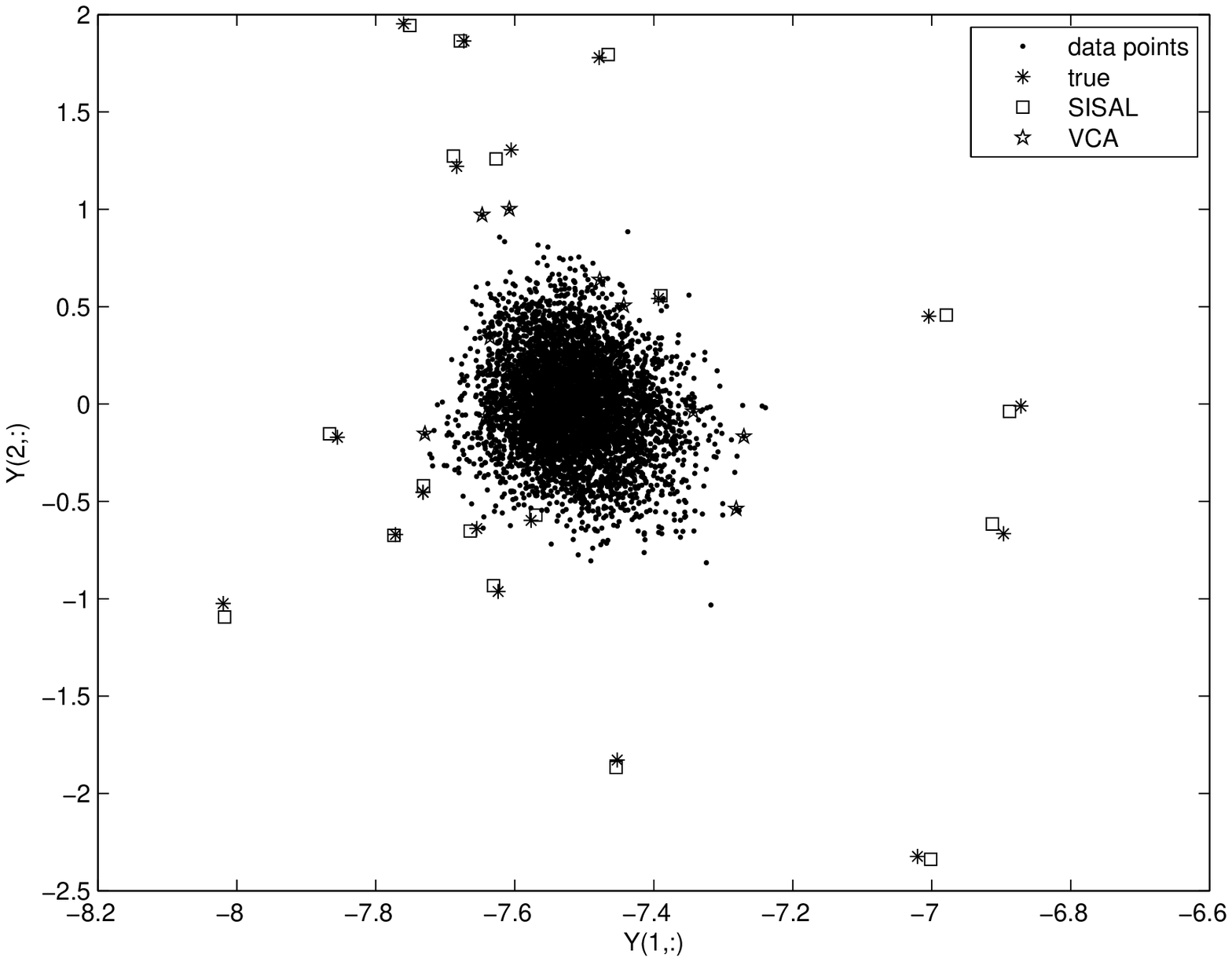}
    }
  \caption{Unmixing results for (a) $p=3$ and (b) $p=20$ number of endmembers for
           SISAL, MVSA, and MVES algorithms. Dots represent spectral vectors; all other
           symbols represent inferred  endmembers by the unmixing algorithms.  The unmixing problem with $n=10000$  spectral vectors and $p=20$ endmembers is far  beyond the reach of MVSA, and MVES.}
  \label{fig:p3_snr0}
\end{figure}

Fig. \ref{fig:p3_snr0} shows a projection on a subspace
of the true endmembers, of the  SISAL, MVSA, and MVES estimates and the
spectral vectors. The data set has size $n = 10000$  and  a number of endmembers
$p=3$, (part a),  and $p = 20$, (part b). In part b,  we just plot SISAL results because
MVES takes  hours for $p\geq 10$, and MVSA  exhausts the memory for $p\geq 12$.

Notice the high quality of SISAL, MVSA, and MVES  estimates  in both
scenarios. This is not the case with VCA,  as it was not conceived for non-pure pixel scenarios.  VCA plays, however, a  valuable rule in the SISAL, MVSA, and MVES  initializations.

\begin{table}
\centering
\caption{Comparison of SISAL, MVSA, and MVES algorithms for different number of endmembers and
sample size $n=10000$. The time is in seconds, the symbol ``*'' means the algorithm ran  out of memory, while $\dagger$
indicates that the algorithm was aborted before converging.
Note the $O(np)$  time complexity of SISAL.}
\vspace{0.1cm}
\begin{tabular}{c|c|c|c|c|c|c}
         \hline
          & \multicolumn{2}{|c|}{SISAL} &  \multicolumn{2}{|c}{MVSA}&  \multicolumn{2}{|c}{MVES} \\
          \cline{2-7}
         p & $\| \bm{\varepsilon}\|_F$ & T & $\| \bm{\varepsilon}\|_F$ & T & $\|\bm{\varepsilon}\|_F$ & T\\\hline
         3 & 0.03  & 3 & 0.02 & 2 & 0.03 & 2 \\
         6 & 0.08  & 4 & 0.10 & 4 & 0.10 & 56 \\
         8 & 0.07 & 10 & 0.18 & 10 & 0.24 & 296  \\
         10 & 0.13 & 10 & 0.25 & 24 & $\dagger$  & $>$1500  \\
         12 & 0.15 & 14 & * & * & $\dagger$  &  $\gg$ 1500 \\
         20 & 0.18 & 16 & * & * & $\dagger$  &  $\gg$ 1500  \\

        \hline
\end{tabular}
\label{tab:compar}
\end{table}

Table \ref{tab:compar}  shows the times in seconds and the  Frobenius norm  $\| \bm{\varepsilon}\|_F$ of the endmember error matrices  $\widehat{
\bf M}-{\bf M}$. The experiments were performed on an PC equipped with a Intel  Core Duo 3GHz CPU and 4 GB of RAM.  The errors are comparable. However, MVES takes much longer and we could not  run it for more than $p=10$. MVSA ran out or memory  for $p>10$. The time SISAL takes is  well approximated by a $O(np)$ bound, what could be inferred from its structure.

\section{conclusions}
\label{sec:conclusion}

SISAL, a new algorithm for hyperspectral unmixing method of minimum volume  class, was introduced.
The unmixing  is achieved by finding  the minimum volume simplex containing the hyperspectral data. This optimization problem was solved by a sequence of variable splitting  augmented Lagrangian optimizations. The algorithm  complexity is  $O(np)$, where $n$ is the number of spectral vectors and $p$ is the number of endmembers what is much faster than the previous state-of-the-art, allowing   to solve problems far beyond the reach of SISAL's competitors.


\begin{center}
{\bf  ACKNOWLEDGMENT}
\end{center}

The author thanks T.-H. Chan, C.-Y. Chi, and W.-K. Ma, for providing the code for the MVES algorithm \cite{Chan_et_al_2009}.

\end{document}